\documentclass{article}

\newtheorem{thm}{Theorem}[section]

\newtheorem{lem}[thm]{Lemma}

\newcommand{\be}{\begin{equation}}
\newcommand{\ee}{\end{equation}}
\newcommand{\ben}{\begin{enumerate}}
\newcommand{\een}{\end{enumerate}}
\newcommand{\beq}{\begin{eqnarray}}
\newcommand{\eeq}{\end{eqnarray}}
\newcommand{\beqn}{\begin{eqnarray*}}
\newcommand{\eeqn}{\end{eqnarray*}}

\newcommand{\pa}{\partial}
\newcommand{\V}{{\rm V}}

\newcommand{\R}{{\rm R}}

\newcommand{\pxi}{ {\pa \over \pa x^i}}

\newcommand{\pyi}{{\pa \over \pa y^i}}

\font\BBb=msbm10 at 10pt
\newcommand{\Bbb}[1]{\mbox{\BBb #1}}

\newcommand{\qed}{\hspace*{\fill}Q.E.D.}  
\title{\Large \bf Funk Metrics and R-Flat Sprays }
\author{Zhongmin Shen}
\date{Revised in June, 2001}

\begin{document}
\maketitle

\begin{abstract}
The well-known Funk metric $F(x,y)$ is projectively flat  with constant flag curvature ${\bf K}=-1/4$ and the Hilbert metric $F_h(x,y):=(F(x,y)+F(x,-y))/2$ is projectively flat with constant curvature ${\bf K}=-1$. These metrics are the special solutions to Hilbert's Fourth Problem. In this paper, we  construct a non-trivial R-flat spray using the Funk metric. It is then an inverse problem in the calculus of variation to find a Finsler metric that induces the R-flat spray. We find an explicit solution to this inverse problem and obtain a non-trivial projectively flat Finsler metric with ${\bf K}=0$.
\end{abstract}
\section{Introduction}
One of fundamental problems in Finsler geometry is to find and study
non-trivial  Finsler metrics of constant (flag) curvature. Several 
non-trivial Finsler metrics of constant curvature have been found. The simplest ones are the Funk metrics and the Hilbert metrics.
The Funk metric $F=F(x,y)$ on a strongly convex domain $\Omega $ in $\R^n$  is defined by
\be
 { x} + {{y}\over F} \in \pa \Omega, \ \ \ \ \ \ y\in T_{ x}\Omega=\R^n.\label{FFunk}
\ee
The Funk metric  is non-reversible,  positively complete and projectively flat with ${\bf K}=-1/4$.
The Hilbert metric on $\Omega$ is obtained from the Funk metric by symmetrization,
\be
F_h (x,y):= {1\over 2} \Big ( F(x, y)+F(x, -y) \Big ), \ \ \ \ \ y\in T_{x}\Omega=\R^n. \label{Hilbert}
\ee
 It is reversible, complete and projectively flat with ${\bf K}=-1$.

Five years ago, R. Bryant   constructed a family of Finsler metrics on  ${\rm S}^n$  of  constant curvature ${\bf K}=1$ (see \cite{Br1}\cite{Br2}\cite{Br3}). The Bryant metrics are non-reversible and projectively flat (geodesics are great circles). 
In \cite{Sh3}, we construct a new family of  pointwise projectively flat Randers metrics  on the unit ball $\Bbb B^n$ with ${\bf K}=-1/4$. They are given by
\be
 F_a ({\bf x}, {\bf y})) = {\sqrt{ |y|^2 - ( |x|^2 |y|^2 - \langle x, y \rangle^2 )}  + \langle x, y \rangle \over 1-|x|^2} 
+ { \langle a, y \rangle \over 1 + \langle a, x \rangle }, \ \ \ \ \ {y}\in T_{x} \R^n =\R^n,\label{Funktype}
\ee
where $a\in \R^n$ is a constant vector with $|a| < 1$, and $|\cdot|$ and $\langle, \rangle $ denote the Euclidean norm and inner
product in $\R^n$, respectively. When $a=0$, the metric $F_{0}$
is just the Funk metric on $\Bbb B^n$.

In \cite{Sh2}, we conjecture that there are
non-trivial positively complete, projectively flat  Finsler metrics of constant curvature ${\bf K}=0$. In this paper, we will prove the   existence of projectively flat  Finsler metrics of  curvature ${\bf K}=0$ by constructing a projectively flat and R-flat spray using the Funk metrics. 

What are sprays?
A spray on a manifold $M$ is a global vector field ${\bf G}$ on $TM$ which is expressed  in a standard local coordinate system  $(x^i, y^i)$ by
\[ {\bf G} = y^i\pxi - 2 G^i(x,y) \pyi,\]
where $G^i(x,y)$ are local $C^{\infty}$ functions on $TM\setminus\{0\}$ satisfying
$G^i(x,\lambda y) = \lambda G^i(x,y)$, $ \lambda >0$.  Every Finsler metric 
induces a spray (see (\ref{Gi})). The notion of Riemann curvature is defined for  sprays \cite{Bw}\cite{Dg}\cite{Ko}\cite{Sh1}.
 A spray  is said to be {\it R-flat} if its Riemann curvature vanishes.
A Finsler metric is of constant curvature ${\bf K}=0$ if and only if its spray is R-flat. 
According to \cite{GrMu}, every isotropic spray is locally induced by 
a Finsler metric. Thus  as long as we find a R-flat spray,  we obtain a Finsler metric with ${\bf K}=0$.

\begin{thm}\label{thm1}
Let $\Omega$ be a strongly convex domain in $\R^n$ and $F$ the 
Funk metric on $\Omega$. Then the following spray is R-flat.
\be
 \tilde{\bf G}:= y^i\pxi -2 F y^i\pyi.\label{Gdelta}
\ee
\end{thm}

From (\ref{Gdelta}), we see that the geodesics of $\tilde{\bf G}$ are straight lines in $\Omega$. 
Thus if a Finsler metric  $\tilde{F}$  induces $\tilde{\bf G}$ on an open subset of $\Omega$, then it is  pointwise projectively flat with ${\bf K}=0$.

\begin{thm}\label{thm3}
Let $F=F(x,y)$ be the Funk metric on a strongly convex domain $\Omega \subset \R^n$. For an arbitrary point $a=(a^i)\in \Omega$, define
a function $\tilde{F}=\tilde{F}(x, y)$ on $T\Omega=\Omega\times \R^n$ by
\be
\tilde{F}: = F (x, y) + F_{x^i}(x, y) (x^i-a^i).\label{solution}
\ee
$\tilde{F}$ is a pointwise projectively flat Finsler metric  with ${\bf K}=0$.
\end{thm}

If we take the Funk metric on $\Bbb B^n$ and $a=0$, then the resulting  Finsler metric $\tilde{F}=\tilde{F}(x,y)$ on $\Bbb B^n$ is given by
\be
\tilde{F} =  {\Big (\sqrt{|y|^2-(|x|^2|y|^2-\langle x, y\rangle^2)} +\langle x, y\rangle\Big )^2 \over (1-|x|^2)^2\sqrt{|y|^2- (|x|^2|y|^2-\langle x, y\rangle^2)}} ,  \label{FFF1}
\ee
where $y\in T_x\Bbb B^n=\R^n$.
By Theorem \ref{thm2}, $\tilde{F}$  is a pointwise projectively flat Finsler metric on $\Bbb B^n$ with ${\bf K}=0$.

Finsler metrics given in (\ref{solution}) are not all Finsler metrics that induce $\tilde{\bf G}$. 
We can express all analytic  Finsler metrics that induces $\tilde{\bf G}$  in a power series.
See Theorem \ref{thm2} below.

\section{Preliminaries}
A  Minkowski norm $\varphi$  on a vector space $\V$ is a
nonnegative function 
with the following properties
\ben
\item[(a)]
 $\varphi$ is positively homogeneous of degree one, i.e.,
\be
\varphi(\lambda {\bf y}) = \lambda \varphi({\bf y}), \ \ \ \ \ \ \lambda >0,\ {\bf y}\in \V.
\ee
\item[(b)] $\varphi$ is $C^{\infty}$ on $\V\setminus\{0\}$
and for any ${\bf y}\in \V\setminus\{0\}$, 
\be
g_{\bf y}({\bf u}, {\bf v}) 
: = {1\over 2} {\pa^2 \over \pa s \pa t} \Big [ \varphi^2({\bf y}+s{\bf u}+t{\bf v})\Big ]_{|s=t=0}, \ \ \ \ \ \ {\bf u}, {\bf v}\in \V,
\ee
is a positive definite symmetric  bilinear form.
\een
A domain $\Omega$ in a vector space $\V$ is said to be {\it strongly convex} if 
there is a point $x_o\in \Omega$ and a Minkowski norm $\varphi$ on $\V$ 
such that $\pa \Omega -\{x_o\} =\varphi^{-1}(1)$. 

Let $M$ be an $n$-manifold. A family of Minkowski norms $F=\{F_p\}_{p\in M}$ in tangent spaces $T_pM$
is called a {\it Finsler metric} on $M$ if it is $C^{\infty}$ on $TM\setminus\{0\}$. Throughout this paper,
Finsler metrics are always positive definite, unless otherwise stated.

\bigskip
\noindent {\bf Funk metrics}: 
Let  $\Omega$ be a strongly convex domain in  $\R^n$. By definition, 
there is a Minkowski norm $\varphi$ on $\R^n$  and a point ${x}_o\in \Omega$ 
such
that $\pa \Omega -\{{x}_o\} = \varphi^{-1}(1)$.
Let $F$ be the Funk metric on
$\Omega$.
For any ${y}\in T_{x}\Omega=\R^n$,  $F= F(x, y)$ is determined by
\be
 \varphi \Big ( { x}-{x}_o + {{y}\over F} \Big ) = 1.\label{phi}
\ee
Differentiating (\ref{phi}) yields a system of PDEs, 
\be
F_{x^i} = FF_{y^i},\ \ \ \ \ \ i=1,\cdots, n.\label{Fxy}
\ee
Equation (\ref{Fxy}) is proved in \cite{Ok}. 
The Funk metric $F$ on the unit ball $\Bbb B^n\subset \R^n$ is given by 
\be
 F = {\sqrt{ |{y}|^2 - ( |{ x}|^2 |{y}|^2 - \langle { x}, {y} \rangle^2 )}  + \langle {x}, {y} \rangle \over 1-|{x}|^2}. \label{Funk}
\ee

\bigskip

Every Finsler metric $F$ on $M$ induces a spray ${\bf G}
= y^i\pxi-2G^i(x,y)\pyi$
by
\be
G^i(x,y) := {1\over 4} g^{il}(x,y) \Big \{ 2 {\pa g_{jl}\over \pa x^k}(x,y) - {\pa g_{jk}\over \pa x^l}(x,y)  \Big \} y^j y^k, \label{Gi}
\ee
where $g_{ij}(x,y):= {1\over 2} [F^2]_{y^iy^j}(x,y)$.

For a vector ${\bf y}=y^i\pxi|_p\in T_pM$, set ${\bf R}_{\bf y}({\bf u}) :=R^i_{\ k} u^k \pxi|_p$,  where ${\bf u}=u^i\pxi|_p$ and $R^i_{\ k}=R^i_{\ k}(x,y)$ are 
given by
\be
 R^i_{\ k} : = 2 {\pa G^i\over \pa x^k}- y^j {\pa^2 G^i\over \pa x^j\pa y^k}  + 2 G^j {\pa^2 G^i\over \pa y^j \pa y^k}-{\pa G^i\over \pa y^j}{\pa G^j\over \pa y^k}.\label{Rik}
\ee
Clearly, 
\be
{\bf R}_{\bf y}({\bf y})=0. \label{Ry=0}
\ee

Assume that  ${\bf G}$ is induced by a Finsler metric $F$, then 
${\bf R}_{\bf y}$ is self-adjoint with respect to $g_{\bf y}$, i.e.,
\be
 g_{\bf y} ( {\bf R}_{\bf y} ({\bf u}),\; {\bf v})=g_{\bf y}({\bf u},\; {\bf R}_{\bf y} ({\bf v}) ).\label{Rsym}
\ee
For a tangent plane $P\subset T_pM$ and a vector ${\bf y}\in P\setminus\{0\}$,
the {\it flag curvature} ${\bf K}(P, {\bf y})$ is defined by
\[ {\bf K} (P, {\bf y}):= { g_{\bf y} ({\bf R}_{\bf y}({\bf u}), {\bf u}) \over g_{\bf y}({\bf y}, {\bf y})g_{\bf y}({\bf u}, {\bf u})-g_{\bf y}({\bf y}, {\bf u})g_{\bf y}({\bf y}, {\bf u})},\]
where ${\bf u}\in P$ such that $P={\rm span}\{ {\bf y}, {\bf u}\}$. By
(\ref{Ry=0}) and (\ref{Rsym}), we see that  ${\bf K}(P, {\bf y})$ is independent of the choice of ${\bf u}\in P$.
Clearly, the flag curvature is a constant,  ${\bf K}= \lambda $ if and only if
\[ {\bf R}_{\bf y}({\bf u}) = \lambda \Big \{ g_{\bf y}({\bf y}, {\bf y})\; {\bf u}-g_{\bf y}({\bf y}, {\bf u})\; {\bf y} \Big \},
\ \ \ \ \ \ {\bf y}, {\bf u }\in T_pM.\]
In particular,
$ {\bf K} =0$ if and only if $ {\bf R}=0$.

\bigskip

\noindent
{\it Proof of Theorem \ref{thm1}}:
Let $F= F(x, y)$ be the Funk metric on a strongly convex domain $\Omega \subset \R^n$ and
\be
\tilde{\bf G}: =y^i\pxi -2  F y^i \pyi.\label{Gdelta*}
\ee
Let
\[
 G^i(x,y) :=  F y^i .
\]
Using (\ref{Fxy}), we obtain
\begin{eqnarray*}
{\pa G^i\over \pa x^k} & = &  F F_{y^k} y^i\\
{\pa G^i\over \pa y^j} & = &  F_{y^j} y^i + F ^i_j\\
y^j {\pa^2 G^i\over \pa x^j \pa y^k} & = & 
 F F_{y^k} y^i + F^2 \delta^i_k\\
G^j {\pa^2 G^i\over \pa y^j \pa y^k} & = &  F F_{y^k} y^i + F^2 \delta^i_k.
\end{eqnarray*}
Plugging them into 
(\ref{Rik}) yields
\[
R^i_{\ k} =0.
\]
Thus $\tilde{\bf G}$ is R-flat.
This proves Theorem \ref{thm1}.
\qed

\bigskip

According to \cite{GrMu}, $\tilde{\bf G}$ can be induced by a Finsler metric. 
In the following sections, we are going to find the solutions 
to this inverse problem in a direct way.

\section{Solving the inverse problem}
We first convert  the above inverse problem to solving  
 a simple system of PDEs. First we need the following 

\begin{lem} {\rm (Rapcs\'{a}k \cite{R})} Let
$\tilde{F}=\tilde{F}(x, y)$ be a Finsler metric on an open subset 
${\cal U}\subset \R^n$. 
$\tilde{F}$ is pointwise projectively flat (i.e., geodesics are straight lines)
if and only if $\tilde{F}$ satisfies
\be
\tilde{F}_{x^k y^l} y^k = \tilde{F}_{x^l}, \ \ \ \ \ \ l=1, \cdots, n.\label{R1}
\ee
In this case, the spray coefficients $G^i$ are in the form $G^i = P y^i$, where
\be
P = {\tilde{F}_{x^k}y^k \over 2 \tilde{F}} .\label{R2}
\ee
\end{lem}
See \cite{Sh1} for details.

\bigskip

Now we are going to prove our key lemma. 
\begin{lem}\label{propPDE}
Let $F$ be the Funk metric on a strongly convex domain $\Omega $ in $\R^n$
and 
\[\tilde{\bf G}= y^i\pxi - 2 F y^i\pyi.\]
Then a Finsler metric 
 $\tilde{F}$ on $\Omega$  induces  $\tilde{\bf G}$ if and only if 
 $\tilde{F}$ satisfies
\be
\tilde{F}_{x^k} =   (F \tilde{F} )_{y^k}, \ \ \ \ \ k=1, \cdots, n.\label{R6}
\ee
\end{lem}
{\it Proof}: Suppose that $\tilde{F}$ induces $\tilde{\bf G}$,
that is, the geodesic coefficients of $\tilde{F}$ are given by
$G^i = F y^i$. 
Since $\tilde{\bf G}$ is pointwise projectively flat, so is $\tilde{F}$. 
By the Rapcs\'{a}k lemma, $\tilde{F}$ satisfies (\ref{R1}) and the geodesic
coefficients $G^i = P y^i$ of $\tilde{F}$ is given by (\ref{R2}). 
Thus
\be
2  F \tilde{F} = \tilde{F}_{x^k} y^k.\label{R3}
\ee
Differentiating (\ref{R3}) with respect to $y^l$  and using (\ref{R1}), we obtain
\[
2  ( F \tilde{F})_{y^l} = \tilde{F}_{x^l} + \tilde{F}_{x^k y^l} y^k
= 2 \tilde{F}_{x^l}.
\]
That is, $\tilde{F}$ satisfies (\ref{R6}).

Conversely, assume that a Finsler metric  $\tilde{F}
= \tilde{F}(x,y)$ on $\Omega$  satisfies
(\ref{R6}). Differentiating (\ref{R6}) 
with respect to $y^l$ and then contracting it with $y^k$ yield
\[ \tilde{F}_{x^ky^l} y^k 
=  (F\tilde{F})_{y^k y^l} y^k
=  (F\tilde{F})_{y^l} =\tilde{F}_{x^l}.\] 
Thus $\tilde{F}$ satisfies (\ref{R1}) and  $\tilde{F}$ is pointwise projectively flat with $G^i = P y^i$ given by (\ref{R2}).
Contracting (\ref{R6}) with $y^l$ yields
\[ \tilde{F}_{x^l}y^l =  (F\tilde{F})_{y^l}y^l = 2  F \tilde{F}.\]
Thus 
\[ P= { \tilde{F}_{x^l}y^l\over 2 \tilde{F}} =  F.\]
Namely, the spray  of $\tilde{F}$ is  just $\tilde{\bf G}$. 
\qed

\bigskip
\noindent{\it Proof of Theorem \ref{thm3}}: 
Let $F=F(x,y)$ denote the Funk metric on a strongly convex domain 
$\Omega\subset \R^n$ and $a \in \Omega^n$. Let
\[ \tilde{F}= F(x, y)+ F_{x^i}(x, y) (x^i -a^i)
= F(x, y) + F(x, y) F_{y^i}(x, y) (x^i-a^i).\]
To show that $\tilde{F}$ induces $\tilde{\bf G}$, it suffices to
verify that $\tilde{F}$ satisfies (\ref{R6}).

By (\ref{Fxy}), we obtain 
\begin{eqnarray*}
\tilde{F}_{x^k} & = & F_{x^k} + F_{x^k} F_{y^i} (x^i-a^i)+ F F_{x^ky^i} (x^i-a^i) +  F F_{y^k}\\
& = &  F F_{y^k} + FF_{y^k} F_{y^i}(x^i-a^i) + F ( FF_{y^k})_{y^i} (x^i-a^i) + F F_{y^k}\\
& = & 2 F F_{y^k} + 2 F_{y^k} F_{y^i} (x^i-a^i) + F^2 F_{y^k y^i} (x^i-a^i)\\
\end{eqnarray*}
and 
\begin{eqnarray*}
(F\tilde{F})_{y^k}
& = & (F^2 + F^2 F_{y^i}(x^i-a^i) )_{y^k}\\
& = & 2 F F_{y^k} + 2 F F_{y^k} F_{y^i} (x^i-a^i) + F^2 F_{y^ky^i} (x^i-a^i).
\end{eqnarray*}
Thus $\tilde{F}$ satisfies that $\tilde{F}_{x^k} = (F\tilde{F})_{y^k}$.
This proves Theorem \ref{thm3}.
\qed

\bigskip
Taking the Funk metric $F$ on the unit ball $\Bbb B^n \subset \Bbb R^n$ 
in (\ref{Funk}) and $a\in \R^n$ with $|a|< 1$, we obtain the following Finsler metric on
$\Bbb B^n$. 
\begin{eqnarray}
\tilde{F} & = &  {F^2\over \sqrt{|y|^2- (|x|^2|y|^2-\langle x, y\rangle^2)}} 
 -{\langle a, y\rangle F+ \langle a, x\rangle F^2  \over \sqrt{|y|^2-(|x|^2|y|^2-\langle x, y\rangle^2)}} \nonumber\\
& = & { (1- \langle a, x\rangle) F^2 - \langle a, y\rangle F \over
\sqrt{|y|^2-(|x|^2|y|^2-\langle x, y\rangle^2)}}. \label{FFFF1}\end{eqnarray}
By Theorem \ref{thm3}, we know that $\tilde{F}$ is pointwise projectively flat with ${\bf K}=0$.

\section{Analytic Finsler Metrics with ${\bf K}=0$}

The system (\ref{R6}) is crucial in determining analytic Finsler metrics 
that induce $\tilde{\bf G}$.

A Finsler metric $\tilde{F}= \tilde{F}(x, y)$ on an open subset 
$\Omega\subset \Bbb R^n$ is said to be analytic at $x_o\in \Omega$
if it can be expressed as a Taylor series around $x_o$ as follows,
\[ \tilde{F}= \sum_{m=0}^{\infty}\sum_{i_1\cdots i_m=1}^n a_{i_1\cdots i_m} (y) (x^{i_1}-x^{i_1}_o)
\cdots (x^{i_m}-x^{i_m}_o),\]
where $a_{i_{1}\cdots i_m}(y)$ are $C^{\infty}$ functions on $\R^n\setminus\{0\}$ satisfying
\[ a_{i_1\cdots i_m}(\lambda y) = \lambda a_{i_1\cdots i_m}(y), \ \ \ \ \ \lambda >0.\]
Thus $a_0 (y)= \tilde{F}(x_o, y)$ is a Minkowski norm on $\R^n$.

\begin{thm}\label{thm2} Let  $\varphi = \varphi(y)$ be a Minkowdki norm  and $\Omega = \{ y\in \R^n \ | \ \varphi(y) < 1 \}$. 
Let $\tilde{\bf G}$ denote the R-flat spray on $\Omega$ defined in (\ref{Gdelta}).  
If 
$\tilde{F}=\tilde{F}(x,y)$ is a Finsler metric  on a neighborhood of the origin $0\in \Omega$ that induces 
$\tilde{\bf G}$, then $\tilde{F}$ 
is given by
\be
\tilde{F}:
= \sum_{m=0}^{\infty}{1\over m!} 
{d^m \over dt^m}\Big [ \varphi^m (y + tx ) \psi (y+tx ) \Big ]|_{t=0}, \label{FFF} 
\ee
where $\psi(y):=\tilde{F}(0, y)$. 
Conversely, for any Minkowski norm $\psi=\psi(y)$ on $\R^n$, 
the function $\tilde{F}$ defined in (\ref{FFF}) induces  $\tilde{\bf G}$, hence it 
 is a pointwise projectively flat Finsler metric 
with ${\bf K}=0$. \end{thm}

Let $\Omega\subset \R^n$ be a strongly convex defined by a Minkowski norm $\varphi=\varphi(y)$ on $\R^n$ ($\pa \Omega:=\varphi^{-1}(1)$). Let $F=F(x,y)$ denote the Funk metric on $\Omega$. From the definition of $F$,
\[ F(0, y) = \varphi(y), \ \ \ \ \ \ y\in \R^n.\]

Suppose that there is a Finsler metric $\tilde{F}=\tilde{F}(x, y)$ on $\Omega$ satisfying
(\ref{R6}), 
\be
\tilde{F}_{x^k} = (F \tilde{F} )_{y^k}. \label{RR6}
\ee
By (\ref{Fxy}) and (\ref{RR6}), we obtain
\begin{eqnarray*}
 \tilde{F}_{x^ix^j} 
& = &  (F\tilde{F})_{y^ix^j}\\
& = &  ( F_{x^j}\tilde{F}+ F\tilde{F}_{x^j})_{y^i}\\
& = & ( F F_{y^j} \tilde{F} + F (F\tilde{F})_{y^j} )_{y^i}\\
& = & ( F^2 \tilde{F})_{y^iy^j} 
\end{eqnarray*}
By induction, we obtain 
\be
 \tilde{F}_{x^{i_1}\cdots x^{i_m} } 
= ( F^m \tilde{F} )_{y^{i_1}\cdots y^{i_m} }. \label{FGGG}
\ee
Let 
\[\psi (y) := \tilde{F} (0, y), \ \ \ \ \ \ y\in \R^n.\]
This gives
\[ \tilde{F}_{x^{i_1}\cdots x^{i_m} } (0, y)
= \Big [ \varphi^m \psi \Big ]_{y^{i_1}\cdots y^{i_m} }(y).\]
Thus, if $\tilde{F}= \tilde{F}(x, y)$ is analytic in $x$ at $x=0$ for  a fixed $y\not=0$, then it must be given by 
\be
\tilde{F} = 
\sum_{m=0}^{\infty} {1\over m!} 
\sum_{i_1\cdots i_m} \Big [ \varphi^m \psi \Big ]_{y^{i_1}\cdots y^{i_m} }(y) x^{i_1}\cdots x^{i_m}.\label{ff}
\ee
We can also express the above power series in the following form
\be
\tilde{F} :=\sum_{m=0}^{\infty} {1\over m!} {d^m\over dt^m}
\Big [ \varphi^m (y+tx)\psi(y+tx) \Big ]|_{t=0}.\label{ff*}
\ee

\bigskip
Let $\tilde{F}= \tilde{F}(x,y)$ be defined by the power series in (\ref{ff*}). Assume that $\tilde{F}$ is convergent 
for $x$ in a neighborhood of $0\in \Omega$.
We claim that $\tilde{F}$ induces the spray 
$\tilde{\bf G}$ in (\ref{Gdelta}). 
 By Lemma \ref{propPDE},
it suffices to  verify that the
function $\tilde{F}$  satisfies (\ref{RR6}).
Differentiating (\ref{ff*}) with respect to $x^k$,
\begin{eqnarray}
\tilde{F}_{x^k} & = & \sum_{m=1}^{\infty} {1\over m!}
\sum_{i_1\cdots i_m=1}^n \Big [ \varphi^m \psi \Big ]_{y^{i_1}\cdots y^{i_m}}(y)\sum_{j=1}^m x^{i_1} \cdots \delta^{i_j}_k \cdots x^{i_m} \nonumber\\
& = & \sum_{m=1}^{\infty} {m \over m!} \sum_{i_1 \cdots i_{m-1}=1}^n
\Big [ \varphi^m \psi \Big ]_{y^k y^{i_1}\cdots y^{i_{m-1}}} (y) x^{i_1}\cdots x^{i_{m-1}}\nonumber\\
& = & \sum_{m=0}^{\infty} {1\over m !} \sum_{i_1\cdots i_m=1}^n
\Big [\varphi^{m+1}\psi\Big ]
_{y^ky^{i_1}\cdots y^{i_m}} (y) x^{i_1}\cdots x^{i_m}.\label{c1}
\end{eqnarray}

On the other hand, it follows from (\ref{Fxy}) that
\[ F_{x^{i_1}\cdots x^{i_m}}(x,y) = {1\over m+1} \Big [ F^{m+1} \Big ]_{y^{i_1}\cdots y^{i_m}} (x, y).\]
This gives
\[ F_{x^{i_1}\cdots x^{i_m}}(0,y) = {1\over m+1} \Big [ \varphi^{m+1} \Big ]_{y^{i_1}\cdots y^{i_m}} (0, y).\]
Thus the Funk metric $F=F(x,y)$ can be expressed by
\[
F
= \sum_{m=0}^{\infty} { 1\over (m+1)!}\sum_{i_1\cdots i_m=1}^n \Big [ \varphi^{m+1} \Big ]_{y^{i_1}
\cdots y^{i_m} }(y) x^{i_{1}}\cdots x^{i_{m} }.
\]
The power series multiplication gives
\[ F\tilde{F} = \sum_{m=0}^{\infty} {1\over (m+1)!} \sum_{i_1\cdots i_m=1}^n 
\Big [ \varphi^{m+1}\psi \Big ]_{y^{i_1}\cdots y^{i_m}} (y) x^{i_1}\cdots x^{i_m}.\]
Differetiating $F\tilde{F}$ with respect to $y^k$, we obtain
\be
(F\tilde{F})_{y^k} 
=\sum_{m=0}^{\infty} {1\over (m+1)!} \sum_{i_1\cdots i_m=1}^n\Big [\varphi^{m+1}\psi\Big ]
_{y^ky^{i_1}\cdots y^{i_m}} (y) x^{i_1}\cdots x^{i_m}. \label{c2}
\ee
From (\ref{c1}) and (\ref{c2}), we conclude that $\tilde{F}$ indeed satisfies (\ref{RR6}). This proves Theorem \ref{thm2}.
\qed

\bigskip
There are infinitely many choices for $\varphi $ and $\psi$. Thus we obtain
infinitely many projectively flat Finsler metrics with ${\bf K}=0$. 
Taking 
\[ \varphi(y) := |y|=: \psi (y), \ \ \ \ \ \ y\in \R^n,\]
we obtain 
\be
\tilde{F}= \sum_{m=0}^{\infty} {1\over m!} 
{d^m \over dt^m} \Big [ |y+tx|^{m+1}
 \Big ]_{|t=0}.\label{ff3}
\ee
This is just the Finsler metric given in (\ref{FFF1}).

\noindent
Zhongmin Shen

\smallskip

\noindent
Department of Mathematical Sciences

\noindent
Indiana Univ.-Purdue Univ. Indianapolis

\noindent IN 46202-3216

\noindent U.S.A.

\noindent
zshen@math.iupui.edu

\end{document}